\documentclass[a4paper]{article}

\usepackage{amsmath,amsthm,amssymb}
\usepackage{graphicx}
\usepackage{epstopdf}
\usepackage{listings}

\usepackage{algorithm,algpseudocodex}

\usepackage{listings}

\theoremstyle{remark}
\newtheorem*{remark}{Remark}

\newcommand{\A}{\mathbf A}
\newcommand{\B}{\mathbf B}
\newcommand{\C}{\mathbf C}
\newcommand{\CC}{\mathbb C}

\newcommand{\dd}{\mathbf d}
\newcommand{\Dat}{\mathbf D_t^\alpha}
\newcommand{\Dx}{\mathbf D_x}
\newcommand{\Ex}{\mathbf E_x}
\newcommand{\Et}{\mathbf E_t}
\newcommand{\Ft}{\mathbf F_t}
\newcommand{\Fx}{\mathbf F_x}

\newcommand{\U}{\mathbf U}
\newcommand{\Uin}{\mathbf U_{inner}}

\newcommand{\R}{\mathbb{R}}

\newcommand{\uu}{\mathbf u}

\newcommand{\cptt}{D_t^{\alpha}}

\newcommand{\cpt}{D_t^{\alpha}}

\newcommand{\cptnum}{D_{t,num}^{\alpha}}

\newcommand{\cptnumstar}{D_{t,num}^{\alpha,*}}

\begin{document}

\title{Numerical approximation of Caputo-type advection-diffusion equations via Sylvester equations}

\author{Francisco de la Hoz\thanks{francisco.delahoz@ehu.eus} \and Peru Muniain\thanks{Corresponding author: peru.muniain@ehu.eus}}

\date{University of the Basque Country UPV/EHU}

\maketitle

\begin{abstract} In this paper, we approximate numerically the solution of Caputo-type advection-diffusion equations of the form $\cptt u(t,x) = a_1(x)u_{xx}(t,x) + a_2(x)u_x(t,x) + a_3u(t,x) + a_4(t,x)$, where $\cpt u$ denotes the Caputo fractional derivative of order $\alpha\in(0,1)$ of $u$ with respect to $t$,  $t\in[0, t_f]$ and the spatial domain can be the whole real line or a closed interval. First, we propose a method of order $3 - \alpha$ to approximate Caputo fractional derivatives, explain how to implement an FFT-based fast convolution to reduce the computational cost, and express the numerical approximation in terms of an operational matrix. Then, we transform a given Caputo-type advection-diffusion equation into a Sylvester equation of the form $\A\U+ \U\B = \C$, and special care is given to the treatment of the boundary conditions, when the spatial domain is a closed interval. Finally, we perform several numerical experiments that illustrate the adequacy of our approach. The implementation has been done in Matlab, and we share and explain in detail the majority of the actual codes that we have used.

\end{abstract}



\section{Introduction}

In this paper, we are interested in solving numerically Caputo-type advection-diffusion equations having the following form:
\begin{equation}
	\label{e:cptu0}
	\begin{cases}
		\cptt u(t,x) = a_1(x)u_{xx}(t,x) + a_2(x)u_x(t,x) + a_3u(t,x) + a_4(t,x), \quad t\in[0, t_f],
		\cr
		u(x, 0) = u_0(x),
	\end{cases}
\end{equation}
where $\cpt u$ denotes the Caputo fractional derivative of order $\alpha\in(0,1)$ of $u$ with respect to $t$ (see, e.g., \cite{Dimitrov2023} and its references); $u_x$ and $u_{xx}$ denote respectively the first and second derivatives of $u$ with respect to $x$; and $a_1(x)$, $a_2(x)$, $a_3(x)$, $a_4(t,x)$ are at least continuous. The spatial domain can be either the whole real line or a closed interval, i.e., $x\in\mathbb R$ or $[x_a, x_b]$, and in the latter case, boundary conditions on $x = x_a$ and $x = x_b$ must be explicitly imposed.

Given a real function $f(t)$, rational derivatives generalize the idea of standard, integer-order derivatives like $f'(t)$ and $f''(t)$. However, unlike their integer-order counterparts, they are not local, i.e., they do not depend on a neighborhood of $t$, but on the whole domain $[0, t_f]$. Even if there are several definitions of rational derivatives, the Caputo fractional derivative of order $\alpha\in(0,1)$, which is the one that we consider in this paper, is given by:
\begin{equation}
	\label{e:cIxafx}
	\cpt f(t) = \frac{1}{\Gamma(1 - \alpha)}\int_0^t\frac{f'(\tau)}{(t - \tau)^{\alpha}}d\tau.
\end{equation}
Observe that, without loss of generality, we are integrating from $\tau = 0$, but it is possible to integrate from other values. Indeed, the point from which we integrate can appear explicitly, by writing  $D_{0,t}^{\alpha}$. Other possible notations may include the use of a capital C (standing for the name Caputo) as a superscript or subscript, e.g., $^C\cpt$, $^CD_{0,t}^{\alpha}$, etc. However, in this paper, since there is no risk of confusion, we stick to the notation $\cpt$ to denote the Caputo fractional derivative with respect to $t$.

In order to better understand \eqref{e:cIxafx}, we can apply it to the monomial $f(t) = t^\beta$, with $\beta\not = 0$ (when $\beta = 0$, we have trivially $\cpt 1 = 0$). Then,
	\begin{align*}
		\cpt t^\beta & = \frac{\beta}{\Gamma(1 - \alpha)}\int_0^t(t - \tau)^{-\alpha}\tau^{\beta - 1}d\tau = \frac{\beta t}{\Gamma(1 - \alpha)}\int_0^1(t - tz)^{-\alpha}(tz)^{\beta - 1}dz
			\cr
		& = \frac{\beta t^{\beta - \alpha}}{\Gamma(1 - \alpha)}\int_0^1(1 - z)^{1 - \alpha - 1}z^{\beta - 1} dz = \frac{\beta t^{\beta - \alpha}}{\Gamma(1 - \alpha)}\operatorname B(1 - \alpha, \beta)
		\cr
		& = \frac{\beta t^{\beta - \alpha}}{\Gamma(1 - \alpha)}\cdot\frac{\Gamma(1 - \alpha)\Gamma(\beta)}{\Gamma(\beta - \alpha + 1)} = \frac{\Gamma(\beta + 1)}{\Gamma(\beta - \alpha + 1)}t^{\beta - \alpha},
	\end{align*}
	where we have performed the change of variable $\tau = tz$, and $\Gamma(\cdot)$ and $\operatorname B(\cdot, \cdot)$ denote Euler's gamma and beta functions, respectively. Therefore, the Caputo fractional derivative is coherent with integer-order derivatives.
	
In order to approximate numerically \eqref{e:cIxafx}  at a family of nodes $0 = t_0 < t_1 < t_2 < \ldots < t_{N} = t_f$, we observe that
\begin{equation*}
\cpt f(t_j) =
\left\{
\begin{aligned}
& 0, & & j = 0,
	\cr
& \frac{1}{\Gamma(1 - \alpha)}\sum_{l = 0}^{j - 1}\int_{t_l}^{t_{l+1}}\frac{f'(\tau)}{(t_j - \tau)^{\alpha}}d\tau, & & 1 \le j \le N,
\end{aligned}
\right.
\end{equation*}
where, with some notational abuse, $\cpt f(t_j) \equiv \cpt f(t)|_{t = t_j}$.  Then, denoting $f_j \equiv f(x_j)$, the simplest option is to interpolate $f(t)$ at each subinterval $[t_l, t_{l+1}]$ by the straight line that joins the points $(t_l, f_l)$ and $(t_l, f_{l+1})$:
\begin{align}
\label{e:cptftfirst}
\cpt f(t_j) & \approx \frac{1}{\Gamma(1 - \alpha)}\sum_{l = 0}^{j - 1}\int_{t_l}^{t_{l+1}}\left[f_l\frac{\tau - t_{l+1}}{t_l - t_{l+1}} + f_{l+1}\frac{\tau - t_{l}}{t_{l+1} - t_{l}}\right]'   (t_j - \tau)^{-\alpha}d\tau
	\cr
& = \frac{1}{\Gamma(1 - \alpha)}\sum_{l = 0}^{j - 1}\frac{f_{l+1} - f_l}{t_{l+1} - t_{l}}\int_{t_l}^{t_{l+1}}   (t_j - \tau)^{-\alpha}d\tau
	\cr
& = \frac{1}{\Gamma(2 - \alpha)}\sum_{l = 0}^{j - 1}\frac{(f_{l+1} - f_l)((t_j - t_l)^{1 - \alpha} - (t_j - t_{l+1})^{1 - \alpha})}{t_{l+1} - t_l}, \quad 1 \le j \le N.
\end{align}
This idea has been generalized in \cite{highorder3} (see also \cite{highorder1,highorder2}), where, after choosing a family of equally spaced nodes $t_j = hj$, with $h = t_f / N$ and $0\le j \le N$, the authors approximate $f(t)$ at each interval $[t_j, t_{j+1}]$ by a polynomial of degree $r$ that interpolates $f(t)$ at the $r + 1$ nodes $\{t_{j-r+1}, t_{j-r+2}, \ldots, t_j, t_{j+1}\}$, when $j - r + 1 \ge 0$, and by a polynomial of degree $j + 1$ that interpolates $f(t)$ at the $j + 2$ nodes $\{t_0, t_1, \ldots, t_j, t_{j+1}\}$, when $j - r + 1 < 0$.  This allows to obtain an approximation of $\cpt f(t)$ of order $r + 1 - \alpha$ (and, hence, \eqref{e:cptftfirst} is of order $2 - \alpha$).

Observe that choosing interpolating nodes not larger than $t_{j+1}$ enables also to solve numerically a PDE like \eqref{e:cptu0}, by approximating $u$ first at $t_1$, then at $t_2$, and so on, until $t_f$. On the other hand, in this paper, we follow a different approach. More precisely, after denoting $u_{ij} \approx u(t_i, x_j)$, where $t_i$ is a time node and $x_j$ is a spatial node, we transform \eqref{e:cptu0} into a matrix equation, namely a Sylvester equation (see, e.g., \cite{bhatia1997}):
\begin{equation}
\label{e:sylvester}
\A\U + \U\B = \C.
\end{equation}
The idea of transforming PDEs into Sylvester equations has being used in, e.g., \cite{delahozvadillo2013a,delahozvadillo2013b}, but, to the best of our knowledge, it has not been applied explicitly to the case of fractional derivatives in time, and allows approximating the solution of \eqref{e:cptu0} simultaneously at all the time instants considered.

Bearing in mind the previous arguments, the structure of this paper is as follows. In Section~\ref{s:cpt3a}, following the ideas in \cite{highorder3}, we generalize \eqref{e:cptftfirst} and develop an approximation of $\cpt f(t)$ of order $3 - \alpha$, for a given function $f(t)$. More precisely, as explained above, we interpolate $f(t)$ at $\{t_{j-1}, t_{j}, t_{j+1}\}$, when $1 \le j \le N$, but, unlike in \cite{highorder3}, when $j = 0$, we interpolate $f(t)$ at the nodes $\{j_0, j_1, j_2\}$, in order to keep the second degree of the interpolating polynomial at all the intervals. In Section~\ref{s:cptFFT}, we express our numerical approximation in terms of discrete convolutions, which enables to apply the FFT to speed up the calculations (hence reducing the computational time), and to consider extremely large amounts of nodes; this idea has been used to approximate numerically the fractional Laplacian \cite{cayama2022,cuesta2024}, but, to the best of our knowledge, it has not been used in the context of Caputo fractional derivatives. In Section~\ref{s:convergence}, we check numerically the order of convergence $3 - \alpha$, and, in Section~\ref{s:Dat}, we express the numerical approximation in terms of an operational matrix $\Dat$.

Section~\ref{s:numerical} is devoted to solving numerically equations of the form of \eqref{e:cptu0}, by transforming them into a Sylvester system equation like \eqref{e:sylvester}. In Section~\ref{s:unbounded}, we consider the case when the spatial domain is $\mathbb R$, and, in Section~\ref{s:bounded}, when it is a closed interval, and we explain in detail how to implement different types of boundary conditions in the latter case. We carry out the numerical experiments in Section~\ref{s:experiments}: one corresponding to an unbounded domain, and two to a bounded one. In fact, the last numerical experiment is taken from \cite[Example 3]{highorder3}, but we are able to obtain smaller errors.

We want to underline that the method developed in this paper is not an ad hoc one, but it could be applied to other approximations of the Caputo fractional derivative, to other domains (e.g., semiinfinite intervals could be easily considered) and to more general equations than \eqref{e:cptu0}, by introducing, e.g., higher order partial derivatives in space. In this regard, the choice of choosing an approximationg of $\cpt f(t)$ whose order of convergence is $3 - \alpha$ is a deliberate one, and obeys to two reasons: the accuracy that it provides is enough from a practical point of view, and the calculations that appear are not too cumbersome and can be shown explicitly, which allows keeping the focus on the transformation of \eqref{e:cptu0} into \eqref{e:sylvester}, which is the central idea of this paper.

We have used Matlab \cite{matlab}, which is particularly well suited for this kind of problems, and have shared and explained in detail the majority of the actual codes, in order to allow the reproduction of the numerical results, and more importantly, to encourage discussion and further improvements. All the simulations have been run in an Apple MacBook Pro (13-inch, 2020, 2.3 GHz Quad-Core Intel Core i7, 32 GB), except for those corresponding to Figure~\ref{f:errorsaN}, which have been executed in a Mainstream A+ Server AS-2024S-TR SUPERMICRO AMD, with 56 cores, 112 threads, 2.75 GHz, and 128 GB of RAM.

\section{A numerical approximation of $\cpt f(t)$ of order $3 - \alpha$}

\label{s:cpt3a}

Given a function $f(t)$, we want to approximate $\cpt f(t)$ at $t_j = jh$, with $h = t_f / N$, $1 \le j \le N$. Reasoning as in \eqref{e:cptftfirst},
\begin{align}
\label{e:cptft}
\cpt f(t_j) & = \frac{1}{\Gamma(1 - \alpha)}\sum_{l = 0}^{j - 1}\int_{t_l}^{t_{l+1}}\frac{f'(\tau)}{(t_j - \tau)^{\alpha}}d\tau \approx \frac{1}{\Gamma(1 - \alpha)}\sum_{l = 0}^{j - 1}\int_{t_l}^{t_{l+1}}\frac{L_l'(\tau)}{(t_j - \tau)^{\alpha}}d\tau,
\end{align}
where $L_0(t)$ is the Lagrange interpolating polynomial that interpolates $f(t)$ at $t_0 = 0$, $t_1 = h$ and $t_2 = 2h$:
\begin{align}
\label{e:L0t}
L_0(t) & = f_0\frac{t - t_{1}}{t_{0} - t_{1}}\frac{t - t_{2}}{t_{0} - t_{2}} +  f_1\frac{t - t_{0}}{t_1 - t_{0}}\frac{t - t_{2}}{t_1 - t_{2}} + f_{2}\frac{t - t_{0}}{t_{2} - t_{0}}\frac{t - t_{1}}{t_{2} - t_{1}}
	\cr
& = f_0\frac{t^2 - 3ht + 2h^2}{2h^2} -  f_1\frac{t^2 - 2ht}{h^2} + f_{2}\frac{t^2 - ht}{2h^2}.
\end{align}
whereas $L_l(t)$, with $1 \le l \le N - 1$, is the Lagrange interpolating polynomial that interpolates $f(t)$ at $t_{l-1} = (l-1)h$, $t_l = lh$ and $t_{l+1} = (l+1)h$:
\begin{align*}
L_l(t)  &= f_{l-1}\frac{t - t_{l}}{t_{l-1} - t_{l}}\frac{t - t_{l+1}}{t_{l-1} - t_{l+1}} +  f_l\frac{t - t_{l-1}}{t_l - t_{l-1}}\frac{t - t_{l+1}}{t_l - t_{l+1}} + f_{l+1}\frac{t - t_{l-1}}{t_{l+1} - t_{l-1}}\frac{t - t_{l}}{t_{l+1} - t_{l}}
	\cr
& = f_{l-1}\frac{t^2 - (2l+1)ht + (l^2+l)h^2}{2h^2}  -  f_l\frac{t^2 - 2lht + (l^2-1)h^2}{h^2} + f_{l+1}\frac{t^2 - (2l-1)ht + (l^2 - l)h^2}{2h^2}.
\end{align*}
We remark that the polynomial $L_0(t)$ corresponding to the first interval $[t_0, t_1]$ has the same (second) degree as the other intervals, unlike in \cite{highorder3}, where it was
\begin{align}
\label{e:L0tstar}
L_0^*(t) & = f_0\frac{t - t_{1}}{t_{0} - t_{1}} +  f_1\frac{t - t_{0}}{t_1 - t_{0}} = f_0\frac{h - t}{h} +  f_1\frac{t}{h}.
\end{align}
In order to compute \eqref{e:cptft}, we bear in mind that
\begin{align*}
\int_{t_l}^{t_{l+1}}\frac{\tau d\tau}{(t_j - \tau)^{\alpha}} & = \int_{t_j - t_{l+1}}^{t_j - t_l}\frac{t_j - \tau}{\tau^{\alpha}}d\tau = \frac{t_j(t_j - t_l)^{1 - \alpha} - t_j(t_j - t_{l+1})^{1 - \alpha}}{1 - \alpha} - \frac{(t_j - t_l)^{2 - \alpha} - (t_j - t_{l+1})^{2 - \alpha}}{2 - \alpha}
	\cr
& = h^{2-\alpha}\frac{j(j - l)^{1 - \alpha} - j(j - l - 1)^{1 - \alpha}}{1 - \alpha} - h^{2-\alpha}\frac{(j - l)^{2 - \alpha} - (j - l - 1)^{2 - \alpha}}{2 - \alpha},
	\cr
\int_{t_l}^{t_{l+1}}\frac{d\tau}{(t_j - \tau)^{\alpha}} & = \frac{(t_j - t_l)^{1 - \alpha} - (t_j - t_{l+1})^{1 - \alpha}}{1 - \alpha} = h^{1 - \alpha}\frac{(j - l)^{1 - \alpha} - (j - l - 1)^{1 - \alpha}}{1 - \alpha}.
\end{align*}
Then,
\begin{align*}
& \int_{t_0}^{t_{1}}\frac{L_0'(\tau)}{(t_j - \tau)^{\alpha}}d\tau = \int_{t_0}^{t_1}\bigg[f_0\frac{2\tau - 3h}{2h^2} -  f_1\frac{2\tau - 2h}{h^2} + f_{2}\frac{2\tau - h}{2h^2}\bigg]\frac{d\tau}{(t_j - \tau)^{\alpha}}
	\cr
& \qquad = h^{-\alpha}\bigg[\frac{ f_2 - 2 f_1 + f_0 }{(1-\alpha)(2-\alpha)} (j^{2-\alpha} - (j-1)^{2-\alpha}) - \frac{f_2 - 4 f_1 + 3f_0}{2-2\alpha} j^{1-\alpha} - \frac{f_2 - f_0}{2(1 - \alpha)}(j-1)^{1-\alpha}\bigg]
\end{align*}
and
\begin{align*}
	\int_{t_l}^{t_{l+1}}\frac{L_l'(\tau)}{(t_j - \tau)^{\alpha}}d\tau & = \int_{t_0}^{t_1}\bigg[f_{l-1}\frac{2\tau - (2l+1)h}{2h^2}  -  f_l\frac{2\tau - 2lh}{h^2} + f_{l+1}\frac{2\tau - (2l-1)h}{2h^2}\bigg]\frac{d\tau}{(t_j - \tau)^{\alpha}}
	\cr
& = h^{-\alpha}\bigg[\frac{ f_{l+1} - 2 f_l + f_{l-1} }{(1-\alpha)(2-\alpha)} ((j-l)^{2-\alpha} - (j - l - 1)^{2-\alpha})
\cr
& \qquad \qquad + \frac{f_{l+1} - f_{l-1}}{2(1 - \alpha)}(j-l)^{1-\alpha} - \frac{3f_{l+1} - 4 f_l + f_{l-1} }{2-2\alpha} (j-l-1)^{1-\alpha}\bigg].
\end{align*}
Introducing these expressions into \eqref{e:cptft}, we get, for $1 \le j \le N$,
\begin{align}
\label{e:cptftj}
	\cpt f(t_j) \approx \cptnum f(t_j) & = \frac{h^{-\alpha}}{\Gamma(2 - \alpha)}\bigg[\frac{ f_2 - 2 f_1 + f_0 }{2-\alpha} (j^{2-\alpha}  - (j-1)^{2-\alpha}) - \frac{f_2 - 4 f_1 + 3f_0}{2} j^{1-\alpha}
	\cr
	& \qquad  - \frac{f_2 - f_0}{2}(j-1)^{1-\alpha} + \sum_{l = 1}^{j-1}\bigg[\frac{ f_{l+1} - 2 f_l + f_{l-1} }{2-\alpha} ((j-l)^{2-\alpha} - (j-l-1)^{2-\alpha})
	\cr
	& \qquad + \frac{f_{l+1} - f_{l-1}}{2}(j-l)^{1-\alpha} - \frac{3f_{l+1} - 4 f_l + f_{l-1} }{2}(j-l-1)^{1-\alpha}\bigg]\bigg],
\end{align}
where the sum is taken as zero, when $j = 1$.

For the sake of comparison, we also give the formula equivalent to \eqref{e:cptftj}, but using \eqref{e:L0tstar} in $[t_0, t_1]$, instead of \eqref{e:L0t}, as is done in \cite{highorder3}. Then,
\begin{equation*}
\int_{t_0}^{t_{1}}\frac{(L_0^*)'(\tau)}{(t_j - \tau)^{\alpha}}d\tau = h^{-\alpha}\frac{f_1 - f_0}{1 - \alpha}(j^{1-\alpha} - (j-1)^{1-\alpha}),
\end{equation*}
and, instead of \eqref{e:cptftj}, we have 
\begin{align}
	\label{e:cptftjstar}
	\cpt f(t_j) \approx \cptnumstar f(t_j) & = \frac{h^{-\alpha}}{\Gamma(2 - \alpha)}\bigg[(f_1 - f_0)(j^{1-\alpha} - (j-1)^{1-\alpha})
	\cr
	& \qquad + \sum_{l = 1}^{j-1}\bigg[\frac{ f_{l+1} - 2 f_l + f_{l-1} }{2-\alpha} ((j-l)^{2-\alpha} - (j-l-1)^{2-\alpha})
	\cr
	& \qquad \qquad + \frac{f_{l+1} - f_{l-1}}{2}(j-l)^{1-\alpha} - \frac{3f_{l+1} - 4 f_l + f_{l-1} }{2}(j-l-1)^{1-\alpha}\bigg]\bigg],
\end{align}
where, the sum is again taken as zero, when $j = 1$. Note that, when $j = 0$, $\cptnum f(t_0) = \cptnumstar f(t_0) = 0$.

In Listing~\ref{code:cpt}, we offer the Matlab program \texttt{CaputoDerivative.m}, which serves to illustrate the implementation of \eqref{e:cptftj} and \eqref{e:cptftjstar}, for $f(t) = e^{2t}$, $\alpha = 0.17$, $t_f = 1.2$ and $N = 100$. The program calculates the maximum errors $\max_{1\le j\le N} |\cptnum e^{2t_j} - \cpt e^{2t_j}|$ and $\max_{1\le j\le N} |\cptnumstar e^{2t_j} - \cpt e^{2t_j}|$, where $\cptnum e^{2t_j}$ and $\cptnumstar e^{2t_j}$ are given respectively by \eqref{e:cptftj} and \eqref{e:cptftjstar}, whereas the expression of $\cpt e^{2t_j}$ can be obtained, e.g., by typing in Matlab
\begin{verbatim}
syms a t tau % $\alpha$, $t$ and $\tau$
assume(0<a&a<1) % Impose that $\alpha\in(0,1)$
f=exp(2*tau); % Define $f$
Daf=int(diff(f)/(t-tau)^a,tau,0,t)/gamma(1-a) % Compute % $D^\alpha f$
\end{verbatim}
which yields
\begin{equation}
\label{e:cpte2t}
\cpt e^{2t} = \frac{2^\alpha e^{2t}(\Gamma(1 - \alpha) - \Gamma(1 - \alpha, 2t))}{\Gamma(1 - \alpha)},
\end{equation}
where $\Gamma(\cdot, \cdot)$ denotes the upper incomplete gamma function:
$$
\Gamma(s, x) = \int_{x}^{\infty}t^{s-1}e^{-t}dt,
$$
and hence, $\Gamma(s) = \Gamma(s, 0)$. Note that, in Listing~\ref{code:cpt}, we have precomputed $j^{1-\alpha}$ and $j^{2-\alpha}$, for $0 \le j \le N$, which has a dramatic impact in the execution time.

\begin{lstlisting}[label=code:cpt, language=Matlab, basicstyle=\footnotesize, caption = {Matlab program \texttt{CaputoDerivative.m}, that tests \eqref{e:cptftj} and \eqref{e:cptftjstar}}]
clear
tic
a=0.17; % $\alpha$
tf=1.2; % $t_f$
N=2^17; % $N$
tj=tf*(0:N)'/N; % $t_j$
h=tf/N; % $h$
f=exp(2*tj); % $f(t_j)$
Dtaf=(2^a/gamma(1-a))*(exp(2*tj).*(gamma(1-a)-igamma(1-a,2*tj))); % $D_t^\alpha f$
toc,tic % Elapsed time needed to generate $D^\alpha f$
Dtafnum=zeros(N+1,1); % $D_{t,num}^\alpha f$
Dtafnumstar=zeros(N+1,1); % $D_{t,num}^{\alpha,*}f$
aux1=(0:N).^(1-a); % Precompute $j^{1-\alpha}$
aux2=(0:N).^(2-a); % Precompute of $j^{2-\alpha}$
for j=2:N
    Dtafnum(j+1)=Dtafnum(j+1)...
        +(aux2(j:-1:2)-aux2(j-1:-1:1))*(f(3:j+1)-2*f(2:j)+f(1:j-1))/(2-a)...
        +.5*aux1(j:-1:2)*(f(3:j+1)-f(1:j-1))...
        -.5*aux1(j-1:-1:1)*(3*f(3:j+1)-4*f(2:j)+f(1:j-1));
end
Dtafnumstar(2:N+1)=Dtafnum(2:N+1)+(f(2)-f(1))*(aux1(2:N+1)-aux1(1:N))';
Dtafnum(2:N+1)=Dtafnum(2:N+1)+(f(3)-2*f(2)+f(1))/(2-a)*(aux2(2:N+1)-aux2(1:N))'...
    -.5*(f(3)-4*f(2)+3*f(1))*aux1(2:N+1)'-.5*(f(3)-f(1))*aux1(1:N)';
Dtafnum=(h^(-a)/gamma(2-a))*Dtafnum; % $D_{t,num}^\alpha f$
Dtafnumstar=(h^(-a)/gamma(2-a))*Dtafnumstar; % $D_{t,num}^{\alpha,*}f$
toc % Elapsed time needed to generate $D_{t,num}^\alpha f$ and $D_{t,num}^{\alpha,*}f$
disp(norm(Dtaf-Dtafnum,inf))
disp(norm(Dtaf-Dtafnumstar,inf))
\end{lstlisting}

\begin{figure}[!htbp]
	\centering
	\includegraphics[width=0.5\textwidth]{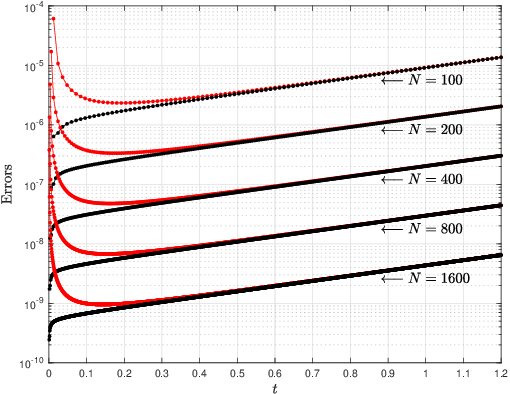}
	\caption{Errors $|\cptnum e^{2t_j} - \cpt e^{2t_j}|$ (in black) and $|\cptnumstar e^{2t_j} - \cpt e^{2t_j}|$ (in red) corresponding respectively to the numerical approximations $\cptnum e^{2t_j}$ (given by \eqref{e:cptftj}) and $\cptnumstar e^{2t_j}$ (given by \eqref{e:cptftjstar}) of $\cpt e^{2t}$, taking $t_f = 1.2$, $\alpha = 0.17$, and $N\in\{100, 200, 400, 800, 1600\}$. }
	\label{f:cptnumoldnew}
\end{figure}

In order to compare the accuracy of \eqref{e:cptftj} and \eqref{e:cptftjstar} in terms of $N$, we plot in Figure~\ref{f:cptnumoldnew}, in semilogarithmic scale, the errors $|\cptnum e^{2t_j} - \cpt e^{2t_j}|$ (in black) and $|\cptnumstar e^{2t_j} - \cpt e^{2t_j}|$ (in red), taking again $t_f = 1.2$, $\alpha = 0.17$, but $N\in\{100, 200, 400, 800, 1600\}$.  As can be seen, the order of convergence of $\cptnum e^{2t_j}$ to $\cpt e^{2t_j}$ and that of $\cptnumstar e^{2t_j}$ to $\cpt e^{2t_j}$ are the same. Moreover, there is virtually no difference between $|\cptnum e^{2t_j} - \cpt e^{2t_j}|$ and $|\cptnumstar e^{2t_j} - \cpt e^{2t_j}|$ for values of $t_j$ close to $t_f = 1.2$, but $\cptnum e^{2t_j}$ is remarkably more accurate than $\cptnum e^{2t_j}$ for small values of $j$. For instance, when $N = 1600$, $|\cptnum e^{2t_1} - \cpt e^{2t_1}| = 1.7425\times10^{-9}$, whereas  $|\cptnumstar e^{2t_1} - \cpt e^{2t_1}| = 1.3460\times10^{-6}$.

In our opinion, the previous results fully justify choosing the interpolating polynomial $L_0(t)$ in \eqref{e:L0t} having the same degree as the other polynomials $L_j(t)$, instead of using $L_0^*(t)$ in \eqref{e:L0tstar}. Therefore, in the rest of the paper, we will use exclusively \eqref{e:cptftj}.

In what regards the order of convergence of \eqref{e:cptftj}, the numerical results reveal that
\begin{align*}
\log_2\left(\frac{|\cptnum(e^{2t_f}, N=100) - \cpt (e^{2t_f}, N=100)|}{|\cptnum(e^{2t_f}, N=200) - \cpt (e^{2t_f}, N=200)|}\right) & = 2.7403,
	\cr
\log_2\left(\frac{|\cptnum(e^{2t_f}, N=200) - \cpt (e^{2t_f}, N=200)|}{|\cptnum(e^{2t_f}, N=400) - \cpt (e^{2t_f}, N=400)|}\right) & = 2.7574,
	\cr
\log_2\left(\frac{|\cptnum(e^{2t_f}, N=400) - \cpt (e^{2t_f}, N=400)|}{|\cptnum(e^{2t_f}, N=800) - \cpt (e^{2t_f}, N=800)|}\right) & = 2.7698,
	\cr
\log_2\left(\frac{|\cptnum(e^{2t_f}, N=800) - \cpt (e^{2t_f}, N=800)|}{|\cptnum(e^{2t_f}, N=1600) - \cpt (e^{2t_f}, N=1600)|}\right) & = 2.7769,
\end{align*}
which is coherent with the theoretical order of convergence $3 - \alpha = 2.83$ predicted by \cite{highorder3}. In the Section~\ref{s:cptFFT}, we make a more careful study on the order of convergence for a large set of values of $\alpha\in(0, 1)$.

\subsection{An FFT-based fast convolution for \eqref{e:cptftj}}

\label{s:cptFFT}

Although the implementation of \eqref{e:cptftj} offered in Listing~\ref{code:cpt} is enough for most purposes, it becomes expensive for larger values of $N$, because its computational cost grows as $\mathcal O(N^2)$, due to the inner sum in \eqref{e:cptftj}. For instance, if $N = 2^{17} = 131072$, the computation of $\cptnum e^{2t_j}$ and $\cptnumstar e^{2t_j}$, for $1 \le j\le N$, requires $72.42$ seconds. Note that the computation of $\cpt e^{2t_j}$, for $1 \le j \le N$, is even more expensive (it needs $201.82$ seconds), because the evaluation of the incomplete $\Gamma$ function \verb|igamma| by Matlab is very costly.

In what follows, we explain how to reduce drastically the computational cost of \eqref{e:cptftj}. 

Given a sequence $a : \mathbb Z\to\CC$ of complex numbers denoted by $\{a_j\}$, and a natural number $P$, $\{a_j\}$ is said to be $P$-periodic if $a_{j + P} =a_j$, for all $j\in\mathbb Z$; note that, in such case, we only need $\{a_0, a_1, \ldots, a_{P-1}\}$ to determine all the values of $a_j$, and also to compute its discrete Fourier transform, which is given by
\begin{equation}
\label{e:FFT}
\hat{a}_{k} = \sum_{j=0}^{P-1}a_je^{-2\pi ijk/P},\quad k\in\mathbb Z.
\end{equation}
It is immediate to check that $\{\hat a_k\}$ is also a $P$-periodic sequence. Moreover, the values $\{a_0, a_1, \ldots, a_{P-1}\}$ and their periodic extension can be recovered by means of the inverse discrete Fourier transform of $\{\hat a_0, \hat a_1, \ldots, \hat a_{P-1}\}$:
\begin{equation}
	\label{e:IFFT}
	a_j = \frac{1}{P}\sum_{k=0}^{P-1}\hat{a}_{k}e^{2\pi ijk/P},\quad j\in\mathbb Z.
\end{equation}
Even if a direct implementation of \eqref{e:FFT} and \eqref{e:IFFT} requires a computational cost of the order of $\mathcal O(P^2)$ operations, \eqref{e:FFT} and \eqref{e:IFFT} can be computed in just $\mathcal O(P\ln(P))$ operations by means of the fast Fourier transform (FFT) and the inverse fast Fourier transform (IFFT), respectively (see \cite{FFT}).

Among the many important applications of the FFT, the one relevant to us in this paper is that it enables to perform efficiently the convolution of two $P$-periodic sequences $\{a_j\}$ and $\{b_j\}$, which is given by
\begin{equation}
\label{e:defconf}
(a\ast b)_{j} = \sum_{l=0}^{P-1}a_{j}b_{l-j},\quad j\in\mathbb Z.
\end{equation}
Then $(a\ast b)_{k} \equiv (b\ast a)_{k}$ is a $P$-periodic sequence satisfying
\begin{equation}
	\label{e:convprop}
	(\widehat{a\ast b})_{k} = \hat{a}_{k}\hat{b}_{k}, \quad k\in\mathbb Z.
\end{equation}
This property, known as the discrete convolution theorem (see, e.g, \cite{cayama2022}), makes possible the computation of \eqref{e:defconf} in just $\mathcal O(P\ln(P))$ operations. On the other hand, if we want to compute an expression of the form 
\begin{equation}
\label{e:convnotperiod}
(a\ast b)_{j} = \sum_{l=0}^{M-1}a_{l}b_{j-l},\quad 0 \le j \le M-1,
\end{equation}
and $\{a_j\}$ and $\{b_j\}$ are not periodic, we cannot apply \eqref{e:convprop} directly, and need to extend them as follows. After observing that we need to know $\{a_0, a_1, \ldots, a_{M-1}\}$ and $\{b_{-M+1}, b_{-M+2}, \ldots, b_{M-1}\}$ to compute \eqref{e:convnotperiod}, we take a period $P$, such that $P \ge 2M-1$, and define
\begin{equation}
	\label{e:tildeab}
	\tilde a_j =
	\begin{cases}
		a_j, & 0 \le j \le M-1,
		\\
		0, & M \le j\le P-1,
	\end{cases}
	\qquad
		\tilde b_j = \begin{cases}
		b_{j}, & 0 \le j\le M-1, \\
		0, & M \le j\le P - M,
		\\
		b_{j-P}, & P - M + 1 \le j\le P - 1.
	\end{cases}
\end{equation}
In this way, we can regard $\{\tilde a_j\}$ and $\{\tilde b_j\}$ as being $P$-periodic. Then, we have immediately that
\begin{equation}
\label{e:bastcP}
(a\ast b)_j \equiv \sum_{l=0}^{M-1}a_{l}b_{j-l} = \sum_{l=0}^{P-1}\tilde a_{l}\tilde b_{j-l} = (\tilde a\ast \tilde b)_{j}, \quad 0 \le j \le M-1,
\end{equation}
where we have used \eqref{e:convprop} to compute $(\tilde a\ast \tilde b)_{j}$. Note that we discard those $(\tilde a\ast \tilde b)_{j}$ for which $M \le j \le P - 1$, because they are not required in \eqref{e:convnotperiod}.

In this paper, we are interested in using this idea to compute efficiently all the values of $\cptnum f(t_j)$ in \eqref{e:cptftj}. As has been mentioned above, the computational cost of \eqref{e:cptftj} is of $\mathcal O(N^2)$ operations, but we can reduce it to only $\mathcal O(N\ln(N))$ operations, by applying the fast convolution theorem to that sum. More precisely, we have to calculate
\begin{align*}
	c_j & = \sum_{l = 1}^{j-1}\bigg[\frac{ f_{l+1} - 2 f_l + f_{l-1} }{2-\alpha} ((j-l)^{2-\alpha} - (j-l-1)^{2-\alpha})
	\cr
	& \qquad \qquad + \frac{f_{l+1} - f_{l-1}}{2}(j-l)^{1-\alpha} - \frac{3f_{l+1} - 4 f_l + f_{l-1} }{2}(j-l-1)^{1-\alpha} \bigg], \quad 2 \le j \le N.
\end{align*}
Then, defining $\tilde c_j \equiv c_{j + 2} \Longleftrightarrow c_j \equiv \tilde c_{j - 2}$, with $0 \le j \le N - 2$,
\begin{align}
\label{e:tildecj}
\tilde c_j & = \sum_{l = 0}^{j}\bigg[\frac{f_{l+2} - 2 f_{l+1} + f_{l} }{2-\alpha} ((j-l+1)^{2-\alpha} - (j-l)^{2-\alpha})
\cr
& \qquad \qquad + \frac{f_{l+2} - f_{l}}{2}(j-l+1)^{1-\alpha} - \frac{3f_{l+2} - 4 f_{l+1} + f_{l}}{2}(j-l)^{1-\alpha}\bigg]
	\cr
& = \sum_{l = 0}^{j}[a_{1,l}b_{1,j-l} + a_{2,l}b_{2,j-l} + a_{3,l}b_{3,j-l}] = \sum_{l = 0}^{N - 2}[a_{1,l}b_{1,j-l} + a_{2,l}b_{2,j-l} + a_{3,l}b_{3,j-l}],
\end{align}
where
\begin{align*}
a_{1,j} & = \frac{f_{j+2} - 2 f_{j+1} + f_{j} }{2-\alpha}, \qquad a_{2,j} = \frac{f_{j+2} - f_{j}}{2}, \qquad a_{3,j} =  -\frac{3f_{j+2} - 4 f_{j+1} + f_{j}}{2},
\cr
b_{1,j} & = ((j+1)^{2-\alpha} - j^{2-\alpha})\chi_{[0,\infty)}(j), \qquad b_{2,j} = (j+1)^{1-\alpha}\chi_{[0,\infty)}(j), \qquad b_{3,j} = j^{1-\alpha}\chi_{[0,\infty)}(j).
\end{align*}
Note that the use of the characteristic function $\chi_{[0,\infty)}(l)$ allows considering $0\le l\le N-2$ in the sum in \eqref{e:tildecj}, instead of $0\le l \le j$, which enables to express \eqref{e:tildecj} as the sum of three convolutions having the form of \eqref{e:convnotperiod}, with $M = N - 1$:
\begin{equation*}
\tilde c_j = (a_1\ast b_1)_{j} + (a_2\ast b_2)_{j} + (a_3\ast b_3)_{j}.
\end{equation*}
Therefore, reasoning as in \eqref{e:tildeab}, we define, for $0 \le j \le N - 2$,
\begin{align*}
\tilde a_{1,j} & = \frac{f_{j+2} - 2 f_{j+1} + f_{j} }{2-\alpha}, \qquad \tilde a_{2,j} = \frac{f_{j+2} - f_{j}}{2}, \qquad \tilde a_{3,j} =  -\frac{3f_{j+2} - 4 f_{j+1} + f_{j}}{2},
\cr
\tilde b_{1,j} & = (j+1)^{2-\alpha} - j^{2-\alpha}, \qquad \tilde b_{2,j} = (j+1)^{1-\alpha}, \qquad \tilde b_{3,j} = j^{1-\alpha},
\end{align*}
and, for $N - 1 \le j \le P - 1$,
$$
\tilde a_{1,j} = \tilde a_{2,j} = \tilde a_{3,j} = \tilde b_{1,j} = \tilde b_{2,j} = \tilde b_{3,j} = 0, \quad 0 \le j \le N - 2,
$$
where $P$ can be any natural number such that $P \ge 2N -3$. In our case, we have chosen $P = 2^{\lceil\log_2(2N -3)\rceil}$, i.e., the smallest power of $2$ satisfying $P \ge 2N -3$, because the FFT and the IFFT are best tuned for powers of $2$. Then,
\begin{equation}
\label{e:tildecjall}
\tilde c_j = (\tilde a_1\ast\tilde  b_1)_{j} + (\tilde a_2\ast\tilde  b_2)_{j} + (\tilde a_3\ast\tilde  b_3)_{j}, \quad 0 \le j \le N - 2,
\end{equation}
and \eqref{e:cptftj} becomes
\begin{align}
\label{e:cptftjnew}
\cptnum f(t_j) & =
\left\{
\begin{aligned}
&  \frac{h^{-\alpha}}{\Gamma(2 - \alpha)}\bigg[\frac{ f_2 - 2 f_1 + f_0 }{2-\alpha} - \frac{f_2 - 4 f_1 + 3f_0}{2} \bigg], & & j = 1,
	\cr
&  \frac{h^{-\alpha}}{\Gamma(2 - \alpha)}\bigg[\frac{ f_2 - 2 f_1 + f_0 }{2-\alpha} (j^{2-\alpha}  - (j-1)^{2-\alpha}) 
	\cr
	& \qquad  - \frac{f_2 - 4 f_1 + 3f_0}{2} j^{1-\alpha} - \frac{f_2 - f_0}{2}(j-1)^{1-\alpha} + \tilde c_{j-2}\bigg], & & 2 \le j \le N.
\end{aligned}
\right.
\end{align}
Note that we need to perform six FFTs, namely those of $\tilde a_1$, $\tilde a_2$, $\tilde a_3$, $\tilde b_1$, $\tilde b_2$, $\tilde b_3$, but only one IFFT, to compute \eqref{e:tildecjall}, which gives a computational cost of $\mathcal O(N\ln(N))$ operations for \eqref{e:tildecjall}, and hence, for \eqref{e:cptftjnew}, too.

In Listing~\ref{code:fastconv}, we offer the Matlab program \texttt{FastConvolution.m}, which tests \eqref{e:cptftjnew} for $f(t) = e^{2t}$, $\alpha = 0.17$, $t_f = 1.2$ and $N = 2^{17} = 131072$. We remark that the elapsed time to compute $\cptnum e^{2t_j}$ by using the fast convolution theorem is only $1.41$ seconds.

\begin{lstlisting}[label=code:fastconv, language=Matlab, basicstyle=\footnotesize, caption = {Matlab program \texttt{FastConvolution.m}, that tests \eqref{e:cptftjnew}}]
clear
tic
a=0.17; % $\alpha$
tf=1.2; % $t_f$
N=2^17; % $N$
tj=tf*(0:N)'/N; % $t_j$
h=tf/N; % $h$
f=exp(2*tj); % $f(t_j)$
Dtaf=(2^a/gamma(1-a))*(exp(2*tj).*(gamma(1-a)-igamma(1-a,2*tj))); % $D_t^\alpha f$
toc,tic % Elapsed time needed to generate $D^\alpha f$
P=2^ceil(log2(2*N-3)); % Length $P$ of $\tilde c_1$ and the other vectors
aux1=(0:N).^(1-a); % Precompute $j^{1-\alpha}$
aux2=(0:N).^(2-a); % Precompute $j^{2-\alpha}$
a1tilde=[(f(3:N+1)-2*f(2:N)+f(1:N-1))/(2-a);zeros(P-N+1,1)]; % $\tilde a_1$
a2tilde=[.5*(f(3:N+1)-f(1:N-1));zeros(P-N+1,1)]; % $\tilde a_2$
a3tilde=[-.5*(3*f(3:N+1)-4*f(2:N)+f(1:N-1));zeros(P-N+1,1)]; % $\tilde a_3$
b1tilde=[(aux2(2:N)-aux2(1:N-1))';zeros(P-N+1,1)]; % $\tilde b_1$
b2tilde=[aux1(2:N)';zeros(P-N+1,1)]; % $\tilde b_2$
b3tilde=[aux1(1:N-1)';zeros(P-N+1,1)]; % $\tilde b_3$
ctilde=ifft(fft(a1tilde).*fft(b1tilde)+fft(a2tilde).*fft(b2tilde)...
    +fft(a3tilde).*fft(b3tilde)); % $\tilde c_1$
Dtafnum=[0;(f(3)-2*f(2)+f(1))/(2-a)-.5*(f(3)-4*f(2)+3*f(1));
    ((f(3)-2*f(2)+f(1))/(2-a))*(aux2(3:N+1)-aux2(2:N))'...
        - .5*(f(3)-4*f(2)+3*f(1))*(2:N)'.^(1-a)...
        - .5*(f(3)-f(1))*aux1(2:N)'+ctilde(1:N-1)];
Dtafnum=(h^(-a)/gamma(2-a))*Dtafnum; % $D_{t,num}^\alpha f$
toc % Elapsed time needed to generate $D_{t,num}^\alpha f$
disp(norm(Dtaf-Dtafnum,inf))
\end{lstlisting}

\subsection{Convergence order of \eqref{e:cptftj} and \eqref{e:cptftjnew}}

\label{s:convergence}

In order to understand how \eqref{e:cptftj} and \eqref{e:cptftjnew} behave as $N$ grows, we have computed $\max_{1\le j\le N} |\cptnum e^{2t_j} - \cpt e^{2t_j}|$, for $\alpha\in\{0.05, 0.1, \ldots, 0.95\}$ and $N\in\{2^1, 2^2, \ldots, 2^{20}\} = \{2, 4, \ldots, 1048576\}$. Since \eqref{e:cptftj} and \eqref{e:cptftjnew} are equivalent, we have used the latter, because its computational cost is much lower. In the literature, in general, it is not that common to consider such large values of $N$, because a given numerical method can be simply unable to deal with them. However, in our opinion, whenever possible, it is interesting to consider values of $N$ as large as possible, to assess correctly a given formula, and have an accurate picture of its behavior.

Since we are taking values of $N$ that are powers of $2$, it is enough to compute once $\cpt e^{2t_j}$ by means of \eqref{e:cpte2t}, for $1 \le j \le N$, and the largest value of $N$, which is $N = 2^{20} = 1048576$ in our case. Then, the values of $\cpt e^{2t_j}$ for smaller values of $N$ are a subset of the values $\cpt e^{2t_j}$ for $N = 2^{20} = 1048576$, so we do not need to recompute $\cpt e^{2t_j}$. We remark that the time required by \eqref{e:cptftjnew} is negligible, even for the whole rank of values of $N$, in comparison to the time necessary to evaluate \eqref{e:cpte2t}.

\begin{figure}[!htbp]
	\centering
	\includegraphics[width=0.5\textwidth]{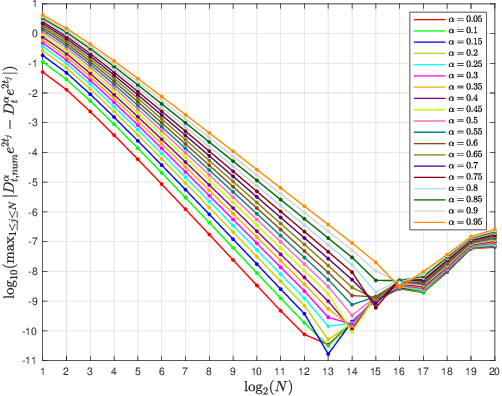}\includegraphics[width=0.5\textwidth]{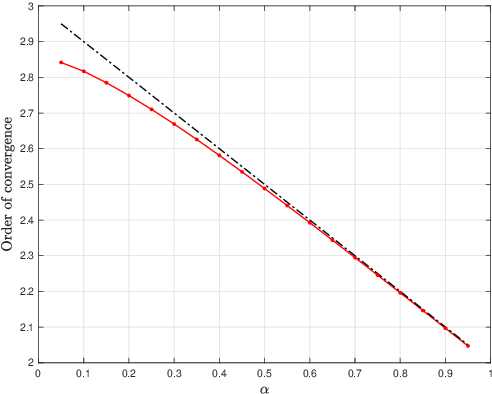}
	\caption{Left: $\log_{10}(\max_{1\le j\le N} |\cptnum e^{2t_j} - \cpt e^{2t_j}|)$ versus $\log_2(N)$, for $\alpha\in\{0.05, 0.1, \ldots, 0.95\}$ and $N\in\{2^1, 2^2, \ldots, 2^{20}\} = \{2, 4, \ldots, 1048576\}$. Right: numerical approximation of the order of convergence (in red) for  $\alpha\in\{0.05, 0.1, \ldots, 0.95\}$ compared to the theoretical values $3 - \alpha$ (dash-dotted black line).}
	\label{f:errorsaN}
\end{figure}

On the left-hand side of Figure~\ref{f:errorsaN}, we have plotted those errors in loglog scale (note that we have conveniently applied binary logarithms  on the $x$-axis, and decimal logarithms on the $y$-axis). The numerical results reveal that all the errors decay as $N$ grows, at least until $N = 2^{13} = 8192$ for all the values of $\alpha$ considered, but, once that the maximum accuracy has been achieved for each value of $\alpha$, the corresponding error starts growing with $N$. Moreover, the maximum possible accuracy is in general better for small values of $\alpha$; for instance, when $\alpha = 0.15$ and $N = 2^{13} = 8192$, $\max_{1\le j\le N} |\cptnum e^{2t_j} - \cpt e^{2t_j}| = 16561\times10^{-11}$, whereas the worst results are obtained when $\alpha = 0.85$, because the lowest error is then $4.9204\times10^{-9}$, which is obtained when $N = 2^{16} = 65536$. In any case, the accuracy diminishes very slowly and, in our opinion, it is completely safe to use \eqref{e:cptftj}, because, even when $N = 2^{20} = 1048576$ (which is such an extremely large value as to be of practical use), the worst error, which corresponds to $\alpha = 0.95$, is $2.6054\times10^{-7}$, so the accuracy is still remarkably high.

The curves on the  left-hand side of Figure~\ref{f:errorsaN} allow also determining approximately the order of convergence in terms of $\alpha$. To do this, we have computed the best fitting line corresponding to $\log_2(N)\in\{8, 9, 10, 11\}$ and their corresponding values $\log_2(\max_{1\le j\le N} |\cptnum e^{2t_j} - \cpt e^{2t_j}|)$ (notice the binary logarithm here). For those values, the Pearson correlation coefficient $\rho$ is very close to $-1$, and the most distant value from $-1$ happens when $\alpha = 0.2$, for which $\rho = -1 + 7.8604\times10^{-7}$. On the right-hand side of   Figure~\ref{f:errorsaN} , we have plotted in red the slope of the regression line for each $\alpha$, together with the theoretical values $3 - \alpha$ forming a black dash-dotted line. The results show that, before the maximum accuracy is reached, the order of convergence is approximately equal to $N^{3 - \alpha}$ (this can be more clearly appreciated when $\alpha$ gets closer to $1$), which is in agreement with the results in \cite{highorder3}.

\subsection{Construction of an operational matrix $\Dat$}

\label{s:Dat}

As we have seen, the use of the fast convolution theorem in \eqref{e:cptftjnew} to speed up the computation of \eqref{e:cptftj} is basically compulsory, if extremely large values of $N$ are to be considered. However, it does not seem straightforward to use  \eqref{e:cptftjnew}, when $f(t)$ is not explicitly known, and, from a practical point of view, it is enough to express \eqref{e:cptftj} as an operational matrix $\Dat$ that multiplies the column vector $(f_0, f_1, \ldots, f_{N})^T$. More precisely, we want to define a matrix $\Dat \equiv [d_{ij}^\alpha]\in\R^{(N+1)\times(N+1)}$, such that \eqref{e:cptftj} becomes
\begin{equation*}
\begin{pmatrix}
\cpt f(t_0) \\ \cpt f(t_1) \\ \vdots \\ \cpt f(t_{N})
\end{pmatrix}
\approx
\begin{pmatrix}
	\cptnum f(t_0) \\ \cptnum f(t_1) \\ \vdots \\ \cptnum f(t_{N})
\end{pmatrix}
=
\Dat \cdot
\begin{pmatrix}
f_0 \\ f_1 \\ \vdots \\ f_{N}
\end{pmatrix}.
\end{equation*}
In order to do it, we  observe in \eqref{e:cptftj} that $f_0$, $f_1$ and $f_2$ intervene in the generation of every $\cptnum f(t_j)$, and that, in general, $\{f_0, f_1, \ldots, f_j\}$ are required to generate $\cptnum f(t_j)$, for $1 \le j \le N$. Therefore, $\Dat$ has the following structure:
$$
\Dat =
\begin{pmatrix}
0 & 0 & 0 & 0 & \ldots & 0
	\cr
d_{21}^\alpha & d_{22}^\alpha & d_{23}^\alpha & 0 &  \ldots & 0
	\cr
d_{31}^\alpha & d_{32}^\alpha & d_{33}^\alpha & 0 & \ldots & 0
	\cr
d_{41}^\alpha & d_{42}^\alpha & d_{43}^\alpha & d_{44}^\alpha & \ldots & 0
	\cr
\vdots & \vdots & \vdots & \vdots & \ddots & \vdots
	\cr
d_{N+1,1}^\alpha & d_{N+1,2}^\alpha & d_{N+1,3}^\alpha & d_{N+1,4}^\alpha & \ldots & d_{N+1,N+1}^\alpha
\end{pmatrix},
$$
i.e., $\Dat$ is almost a lower triangular matrix, except for the entry $d_{23}^\alpha \not=0$, and the nonzero entries of $\Dat$ can be computed in a straightforward way from \eqref{e:cptftj}, by means of Algorithm~\ref{alg:Da}. Note that the indices of a matrix start at $1$, and not at $0$, and that we have omitted the subscript $t$ in the entries $d_{ij}^{\alpha}$, in order not to burden the notation.

\begin{algorithm}
	\caption{Computation  of $\Dat$ \label{alg:Da}}
	\begin{algorithmic}
		\State $\Dat\gets$ array of zeros with size $(N+1)\times(N+1)$
		
		\For{$j\gets 2$ \textbf{to} $N + 1$}

		\State$d_{j1}^\alpha\gets \frac{1}{2-\alpha} ((j-1)^{2-\alpha}  - (j-2)^{2-\alpha}) - \frac{3}{2} (j-1)^{1-\alpha} + \frac{1}{2}(j-2)^{1-\alpha} $

		\State$d_{j2}^\alpha\gets - \frac{2}{2-\alpha} ((j-1)^{2-\alpha}  - (j-2)^{2-\alpha}) + 2 (j-1)^{1-\alpha}$

		\State$d_{j3}^\alpha\gets \frac{1}{2-\alpha} ((j-1)^{2-\alpha}  - (j-2)^{2-\alpha}) - \frac{1}{2} (j-1)^{1-\alpha} - \frac{1}{2}(j-2)^{1-\alpha} $
		
		\EndFor

		\For{$j\gets 3$ \textbf{to} $N + 1$}
		\For{$l\gets 1$ \textbf{to} $j-2$}

	\State 	$d_{jl}^\alpha\gets d_{jl}^\alpha + \frac{1}{2-\alpha} ((j-l-1)^{2-\alpha} - (j-l-2)^{2-\alpha})  - \frac{1}{2}(j-l-1)^{1-\alpha} - \frac{1}{2}(j-l-2)^{1-\alpha}$
	
	\State 	$d_{j,l+1}^\alpha\gets d_{j,l+1}^\alpha - \frac{2}{2-\alpha} ((j-l-1)^{2-\alpha} - (j-l-2)^{2-\alpha}) + 2(j-l-2)^{1-\alpha}$

	\State 	$d_{j,l+2}^\alpha\gets d_{j,l+2}^\alpha + \frac{1}{2-\alpha} ((j-l-1)^{2-\alpha} - (j-l-2)^{2-\alpha}) + \frac{1}{2}(j-l-1)^{1-\alpha} - \frac{3}{2}(j-l-2)^{1-\alpha}$

		\EndFor
		\EndFor

	\State $\Dat\gets\frac{(t_f/N)^{-\alpha}}{\Gamma(2-\alpha)}\Dat$
	
\end{algorithmic}
\end{algorithm}

On the other hand, it is straightforward to implement Algorithm~\ref{alg:Da} in Matlab, which is done in Listing~\ref{code:Da}.

\begin{lstlisting}[label=code:Da, language=Matlab, basicstyle=\footnotesize, caption = {Matlab function \texttt{GenerateDa.m}, that creates $\Dat$ from \eqref{e:cptftj}}]
function Dta=GenerateDa(N,tf,a)
Dta=zeros(N+1); % Preallocate $\mathbf D_t^\alpha$
aux1=(0:N).^(1-a); % Precompute $j^{1-\alpha}$
aux2=(0:N).^(2-a); % Precompute $j^{2-\alpha}$
Dta(2:N+1,1:3)=[(aux2(2:N+1)-aux2(1:N))/(2-a)-1.5*aux1(2:N+1)+.5*aux1(1:N);
    -(2/(2-a))*(aux2(2:N+1)-aux2(1:N))+2*aux1(2:N+1);
    (aux2(2:N+1)-aux2(1:N))/(2-a)-.5*aux1(2:N+1)-.5*aux1(1:N)]';
for j=3:N+1
    Dta(j,1:j-2)=Dta(j,1:j-2)+(aux2(j-1:-1:2)-aux2(j-2:-1:1))/(2-a)...
        -.5*aux1(j-1:-1:2)-.5*aux1(j-2:-1:1);
    Dta(j,2:j-1)=Dta(j,2:j-1)-(2/(2-a))*(aux2(j-1:-1:2)-aux2(j-2:-1:1))...
        +2*aux1(j-2:-1:1);
    Dta(j,3:j)=Dta(j,3:j)+(aux2(j-1:-1:2)-aux2(j-2:-1:1))/(2-a)...
        +.5*aux1(j-1:-1:2)-1.5*aux1(j-2:-1:1);
end
Dta=((tf/N)^(-a)/gamma(2-a))*Dta;
\end{lstlisting}

\section{Numerical approximation of PDEs containing a Caputo de\-riv\-a\-tive in time}

\label{s:numerical}

As said in the introduction, our goal is to approximate numerically the solutions of Caputo-type advection-diffusion equations having the form of \eqref{e:cptu0}:
\begin{equation}
	\label{e:cpttulinear}
	\begin{cases}
		\cptt u(t,x) = a_1(x)u_{xx}(t,x) + a_2(x)u_x(t,x) + a_3u(t,x) + a_4(t,x), \quad t\in[0, t_f],
		\cr
		u(x, 0) = u_0(x),
	\end{cases}
\end{equation}
where $a_1(x)$, $a_2(x)$, $a_3(x)$, $a_4(t,x)$ are at least continuous, and $x\in\mathbb R$ or $[x_a, x_b]$. Note that the ideas in this section can be easily extended to other types of spatial domains like semi-infinite ones, to PDEs containing higher order derivatives with respect to $x$, etc.

The idea is to represent the approximated solution $u(t,x)$ in matrix form, i.e., create a matrix $\U \equiv (u_{ij})\in\CC^{(N_t+1)\times(N_x+1)}$ or $\U \equiv (u_{ij})\in\CC^{(N_t+1)\times N_x}$, such that $u_{ij} \approx u(t_i, x_j)$, where $0 \le i \le N_t$, and $0 \le j \le N_x$ or $0 \le j \le N_x-1$, discretize the partial derivatives $u_x$ and $u_{xx}$ by means of differentiation matrices, impose the initial condition $u(x, 0) = u_0(x)$ and the boundary conditions, if any, and transform \eqref{e:cpttulinear} into a Sylvester equation. In this way, we are able to obtain simultaneously the approximate solution of $\U$ for all the values of $t$.

Note that the variable $t$ changes along the rows of $\U$, whereas the variable $x$ changes along the columns of $\U$.  Therefore, a differentiation matrix must be transposed, in order to differentiate $\U$ with respect to $x$. We denote such transposed matrix as $\Dx$; then, $u_x \approx \U\Dx$, etc.

\subsection{Case with $x\in\mathbb R$}

\label{s:unbounded}

We explain first the case with $x\in\mathbb R$, which is the easiest to implement. In order to discretize $u_x$ and $u_{xx}$, we follow a pseudospectral approach and use Hermite differentiation matrices, although other types of differentiation matrices can be considered with little modification. Those matrices, which, once transposed, are denoted by $\Dx$ and $\Dx^2$, can be generated simultaneously in Matlab by, e.g., the function \verb|herdif| (see \cite{WeidemanReddy2000}), and are based on the Hermite functions:
$$
\psi_n(x) = \frac{e^{-x^2/2}}{\pi^{1/4}\sqrt{2^nn!}}H_n(x),
$$
where $H_n(x)$ are the Hermite polynomials:
$$
H_n(x) = \frac{(-1)^n}{2^n}e^{x^2}\frac{d^n}{dx^n}e^{-x^2},
$$
and the nodes $x_j$, with $0 \le j \le N_x-1$, are the roots of the polynomial $H_n(x)$ with the largest index, multiplied by a scale factor $b$, which is precisely the third parameter of \verb|herdif|.

Let us store the entries of $a_1(x_j)$, $a_2(x_j)$, $a_3(x_j)$ in the diagonals of the diagonal matrices $\A_1$, $\A_2$, $\A_3$, respectively, and $a_4(t_i, x_j)$  in $\A_4$, by doing $a_{4,ij} = a_4(t_i, x_j)$. Then, the right-hand side of \eqref{e:cpttulinear} becomes $a_1(x)u_{xx}(t,x) + a_2(x)u_x(t,x) + a_3(x)u(t,x) + a_4(t,x) \approx \U\B_x + \A_4$, where $\B_x \equiv \Dx^2\A_1 + \Dx\A_2 + \A_3$. On the other hand, the numerical approximation of the left-hand side of \eqref{e:cpttulinear} is given by  $\Dat u \approx \Dat\U$, where $\Dat$ has been constructed in Section~\ref{s:Dat}. At this point, we remark that one advantage of not truncating the spatial domain by considering the whole real line is that it is not necessary to enforce explicitly the boundary conditions at $x = \pm\infty$.  Therefore, the numerical approximation of \eqref{e:cpttulinear} satisfies
\begin{equation}
\label{e:DatUBxA4}
\Dat\U = \U\B_x + \A_4.
\end{equation}
In what regards the initial condition $u_0(x)$, we store the entries $u_0(x_j)$, for $0 \le j \le N_x-1$, in the first row of a matrix $\Ft\in\CC^{(N_t+1)\times N_x}$ whose other elements are zero, and also define a matrix $\Et\in\CC^{(N_t+1)\times N_t}$, whose first row is the all-zero vector of length $N_t$, and the other rows form the identity matrix of order $N_t$:
$$
\Et =
\left(
\begin{array}{ccc}
	0 & \ldots & 0
	\cr
	\hline
	1 & & 0
	\cr
	& \ddots &
	\cr
	0 & & 1
\end{array}
\right),
\qquad
\Ft =
\left(
\begin{array}{ccc}
	u_0(x_0) & \ldots & u_0(x_{N_x + 1})
	\cr
	\hline
	0 & \ldots & 0
	\cr
	\vdots & \ddots & \vdots
	\cr
	0 & \ldots & 0
\end{array}
\right).
$$
Then, we decompose $\U$ as
\begin{equation}
\label{e:Uboundaryt}
\U = \Et\Uin+ \Ft,
\end{equation}
where $\Uin\in\CC^{N_t\times N_x}$ denotes the matrix of $\U$ without its first row, but keeping all its columns:
\begin{equation*}
\Uin =
\begin{pmatrix}
u_{21} & \ldots & u_{2,N_x+1}
	\cr
\vdots & \ddots & \vdots
	\cr
u_{N_t+1,1} & \ldots & u_{N_t+1,N_x+1}
\end{pmatrix}.
\end{equation*}
Introducing \eqref{e:Uboundaryt} into \eqref{e:DatUBxA4},
$$
\Dat(\Et\Uin+ \Ft) = (\Et\Uin+ \Ft)\B_x + \A_4.
$$
Expanding the last expression, and left-multiplying it by $\Et^T$:
$$
\Et^T\Dat\Et\Uin+ \Et^T\Dat\Ft = \Et^T\Et\Uin\B_x + \Et^T\Ft\B_x + \Et^T\A_4,
$$
but $\Et^T\Et$ is the identity matrix of order $N_t$, and $\Et^T\Ft$ is the all-zero matrix. Hence,
$$
\Et^T\Dat\Et\Uin -  \Uin\B_x = -\Et^T\Dat\Ft + \Et^T\A_4.
$$
Finally, denoting
\begin{align*}
\A & \equiv \Et^T\Dat\Et,
	\cr
\B & \equiv - \B_x,
	\cr
\C & \equiv  -\Et^T\Dat\Ft + \Et^T\A_4,
\end{align*}
we arrive at the following Sylvester equation:
\begin{equation}
\label{e:sylv1}
\A\Uin + \Uin\B = \C.
\end{equation}
Once that $\Uin$ has been obtained, $\U$ follows from \eqref{e:Uboundaryt}.

\begin{remark} Sylvester equations like \eqref{e:sylv1} have been extensively studied (see, e.g., \cite{bhatia1997} and its references). Recall that \eqref{e:sylv1} has a unique solution $\Uin$, if and only if the sum of any eigenvalue of $\A$ and any eigenvalue of $\B$ is never equal to zero. Giving a universal proof that the Sylvester equations appearing in this paper satisfy this property for every single case would be difficult and lies beyond the scope of this paper, but there is strong numerical evidence that this indeed happens. On the other hand, it is possible to use the  Bartels-Stewart algorithm \cite{bartels1972} to solve \eqref{e:sylv1}, which requires to obtain the Schur decompositions of $\A$ and $\B$ (observe that this is particularly easy in the case of $\A=\Et^T\Dat\Et$, which is almost triangular, except for the entry $a_{12}\not=0$.) However, from a practical point of view, we have found that the Matlab command \verb|lyap| (which is invoked by typing \verb|lyap(A,B,-C)|) perfectly matches our needs, and we have used it systematically, instead of implementing an algorithm of our own. 
\end{remark}

\subsection{Case with $x\in[x_a, x_b]$}

\label{s:bounded}
	
In this case, we have used again a pseudospectral approach; more precisely, we have considered Chebyshev differentiation matrices \cite{trefethen2000} to discretize $u_x$ and $u_{xx}$, but it is perfectly possible to consider differentiation matrices based on finite differences or on any other technique.

The Chebyshev differentiation matrices corresponding to the interval $[-1, 1]$ can be generated by the function \verb|cheb.m| of \cite{trefethen2000}. They are based on the Chebyshev polynomials $T_n(x) = \cos(n\arccos(x))$, and the nodes $x_j$ are the roots of the polynomial $T_n(x)$ with the largest index, i.e., $x_j = \cos(\pi j/N_x)$, for $0 \le j \le N_x$. In our case, since $x\in[x_a, x_b]$, we have to multiply by $2 / (x_b - x_a)$ the matrix returned by \verb|cheb.m|, which, once transposed, we denote as $\Dx$. Then, the corresponding nodes are accordingly $x_j = x_a + (x_b - x_a)\cos(\pi j/N_x)/2$, for $0 \le j \le N_x$.

Using the new definition of $\Dx$ and following the same steps as in the case with $x\in\R$, we arrive again at \eqref{e:DatUBxA4}, but we have to impose now the boundary conditions:
\begin{equation}
\label{e:mixuatubt}
\begin{cases}
c_au(x_a,t) + d_au_x(x_a,t) = u_a(t), \cr c_bu(x_b, t) + d_bu_x(x_b, t) = u_b(t),
\end{cases}
\end{equation}
where $|c_a| + |d_a| > 0$ and $|c_b| + |d_b| > 0$. Therefore, it is possible to consider Dirichlet boundary conditions, Neumann boundary conditions or Robin boundary conditions. Moreover, the equivalent of \eqref{e:mixuatubt} in matrix form is given by
\begin{equation}
\label{e:mixuatubtmatrix}
\begin{cases}
c_a
\uu_{N_x+1}
+
d_a
\U
\dd_{N_x+1}
=
\uu_a,
\cr
c_b
\uu_1
+
d_b
\U
\dd_1
=
\uu_b,
\end{cases}
\end{equation}
where $\uu_1\in\CC^{N_t+1}$ and $\uu_{N_t+1}\in\CC^{N_t+1}$ are respectively the first and last columns of $\Dx$, and $\dd_1\in\CC^{N_x+1}$ and $\dd_{N_x+1}\in\CC^{N_x+1}$ are respectively the first and last columns of $\Dx$:
$$
\uu_1 =
\begin{pmatrix}
	u_{11} \cr \vdots \cr u_{N_t+1,1}
\end{pmatrix},
\qquad
\uu_{N_x+1} =
\begin{pmatrix}
u_{1,N_x+1} \cr \vdots \cr u_{N_t+1,N_x+1}
\end{pmatrix},
\qquad
\dd_1 = 
\begin{pmatrix}
d_{11} \cr \vdots \cr d_{N_x+1,1}
\end{pmatrix},
\qquad
\dd_{N_x+1} = 
\begin{pmatrix}
d_{1,N_x+1} \cr \vdots \cr d_{N_x+1,N_x+1}
\end{pmatrix},
$$
and the right-hand sides $\uu_a\in\CC^{N_t+1}$ and $\uu_b\in\CC^{N_t+1}$ of \eqref{e:mixuatubtmatrix} are given by
\begin{equation*}
\uu_a =
\begin{pmatrix}
	u_a(t_0)
	\\
	\vdots
	\\
	u_a(t_{N_t})
\end{pmatrix},
\qquad
\uu_b =
\begin{pmatrix}
	u_b(t_0)
	\\
	\vdots
	\\
	u_b(t_{N_t})
\end{pmatrix}.
\end{equation*}
Observe that $\uu_{N_x+1}$ and $\dd_{N_x+1}$ correspond to $\uu_a$, and that $\uu_1$ and $\dd_1$ correspond to $\uu_b$, due to the fact that the Chebyshev nodes are given in reverse order, i.e., $x_a = x_{N_x}$ and $x_b = x_0$. Therefore, should we want to use, e.g., other types of nodes that are given in growing order, it would be enough to swap the subscripts $_a$ and $_b$ in \eqref{e:mixuatubtmatrix}, and we would have
\begin{equation*}
	\begin{cases}
		c_a
		\uu_1
		+
		d_a
		\U
		\dd_1
		=
		\uu_a,
	\cr
		c_b
		\uu_{N_x+1}
		+
		d_b
		\U
		\dd_{N_x+1}
		=
		\uu_b.
	\end{cases}
\end{equation*}
Coming back to \eqref{e:mixuatubtmatrix}, we define $c_{11} = d_ad_{1,N_x+1}$, $c_{12} = c_a + d_ad_{N_x+1,N_x+1}$, $c_{21} = c_b + d_bd_{11}$, $c_{22} = d_bd_{N_x+1,1}$; then, \eqref{e:mixuatubtmatrix} becomes
\begin{equation}
\label{e:mixuatubtmatrixnew}
\begin{cases}
c_{11}
\uu_1
+
c_{12}
\uu_{N_x+1}
=
\uu_a
-
d_a
\U(1:N_t+1,2:N_x)
\dd_{N_x+1}(2:N_x),
\cr
c_{21}
\uu_1
+
c_{22}
\uu_{N_x+1}
=
\uu_b
-
d_b
\U(1:N_t+1,2:N_x)\dd_1(2:N_x),
\end{cases}
\end{equation}
where, using Matlab's syntax, $\U(1:N_t+1,2:N_x)\in\CC^{(N_t+1)\times(N_x-1)}$ denotes $\U$ without its first and last columns, and $\dd_1(2:N_x)\in\CC^{N_x-1}$ and $\dd_{N_x+1}(2:N_x)\in\CC^{N_x-1}$ denote $\dd_1$ and $\dd_{N_x+1}$ without its first and last entries, respectively. Solving \eqref{e:mixuatubtmatrixnew} for $\uu_1$ and $\uu_{N_x+1}$, we get
\begin{align*}
\uu_1 & = \frac{c_{22}\uu_a - c_{12}\uu_b + \U(1:N_t+1,2:N_x)(-d_ac_{22}\dd_{N_x+1}(2:Nx) + d_bc_{12}\dd_1(2:Nx))}{c_{11}c_{22} - c_{12}c_{21}},
	\cr
\uu_{N_x+1} & = \frac{-c_{21}\uu_a + c_{11}\uu_b + \U(1:N_t+1,2:N_x)(d_ac_{21}\dd_{N_x+1}(2:Nx) - d_bc_{11}\dd_1(2:Nx))}{c_{11}c_{22} - c_{12}c_{21}}.
\end{align*}
Bearing in mind the previous arguments, we define,
\begin{equation}
\label{e:Ex}
\Ex =
\left(
\begin{array}{c|ccc|c}
\dfrac{-d_ac_{22}d_{N_x+1,2} + d_bc_{12}d_{12}}{c_{11}c_{22} - c_{12}c_{21}} & 1 & & 0 & \dfrac{d_ac_{21}d_{N_x+1,2} - d_bc_{11}d_{12}}{c_{11}c_{22} - c_{12}c_{21}}
	\\
	\vdots & & \ddots & & \vdots
	\\
\dfrac{-d_ac_{22}d_{N_x+1,N_x} + d_bc_{12}d_{1,N_x}}{c_{11}c_{22} - c_{12}c_{21}} & 0 & & 1 & \dfrac{-d_ac_{22}d_{N_x+1,N_x} - d_bc_{11}d_{1,N_x}}{c_{11}c_{22} - c_{12}c_{21}}
\end{array}
\right)
\end{equation}
and
\begin{equation}
\label{e:Fx}
\Fx =
\left(
\begin{array}{c|ccc|c}
\dfrac{c_{22}u_a(t_0) - c_{12}u_b(t_0)}{c_{11}c_{22} - c_{12}c_{21}} & 0 & \ldots & 0 & \dfrac{-c_{21}u_a(t_0) + c_{11}u_b(t_0)}{c_{11}c_{22} - c_{12}c_{21}}
	\\
	\vdots & \vdots & \ddots & \vdots & \vdots
	\\
\dfrac{c_{22}u_a(t_{N_t}) - c_{12}u_b(t_{N_t})}{c_{11}c_{22} - c_{12}c_{21}} & 0 & \ldots & 0 & \dfrac{-c_{21}u_a(t_{N_t}) + c_{11}u_b(t_{N_t})}{c_{11}c_{22} - c_{12}c_{21}}
\end{array}
\right),
\end{equation}
which enables us to decompose $\U$ as
\begin{equation*}
\U = \U(1:N_t+1,2:N_x)\Ex+ \Fx.
\end{equation*}
Therefore, removing the first row:
\begin{equation}
	\label{e:Uboundaryx}
	\U(2:N_t+1, 1:N_x+1) = \U(2:N_t+1, 2:N_x)\Ex+ \Fx(2:N_t+1, 1:N_x+1).
\end{equation}
Denoting now $\Uin\equiv \U(2:N_t+1, 2:N_x)$, and inserting \eqref{e:Uboundaryx} into \eqref{e:Uboundaryt},
\begin{equation}
	\label{e:Ugeneralboundary}
	\U = \Et\Uin\Ex+ \Et\Fx(2:N_t+1, 1:N_x+1)+ \Ft.
\end{equation}
Inserting \eqref{e:Ugeneralboundary} into \eqref{e:DatUBxA4}:
\begin{align*}
	& \Dat(\Et\Uin\Ex+ \Et\Fx(2:N_t+1, 1:N_x+1)+ \Ft)
	\cr
	& \qquad = (\Et\Uin\Ex+ \Et\Fx(2:N_t+1, 1:N_x+1)+ \Ft)\B_x + \A_4.
\end{align*}
Expanding the last expression, left-multiplying it by $\Et^T$, right-multiplying it by $\Ex^*(\Ex\Ex^*)^{-1}$, where $\Ex^*$ denotes the conjugate transpose of $\Ex$, because it can be complex, and bearing in mind that $\Et^T\Et$ is the identity matrix of order $N_t$ and that $\Et^T\Ft$ is the all-zero matrix,
\begin{align*}
	& \Et^T\Dat\Et\Uin - \Uin\Ex\B_x\Ex^*(\Ex\Ex^*)^{-1} = \big[-\Et^T\Dat\Et\Fx(2:N_t+1, 1:N_x+1)
	\cr
	& \qquad - \Et^T\Dat\Ft +  \Fx(2:N_t+1, 1:N_x+1)\B_x + \Et^T\A_4\big]\Ex^*(\Ex\Ex^*)^{-1}.
\end{align*}
On the other hand, $\Ex^*(\Ex\Ex^*)^{-1}$ is precisely the Moore-Penrose inverse $\Ex^\dagger$ of $\Ex$, because $\Ex$ is a full-rank matrix having more rows than columns. Therefore, we write $\Ex^\dagger$ rather than $\Ex^*(\Ex\Ex^*)^{-1}$, which is advantageous from a numerical point of view. Then, denoting
\begin{align*}
		\A & \equiv \Et^T\Dat\Et,
		\cr
		\B & \equiv -\Ex\B_x\Ex^\dagger,
		\cr
		\C & \equiv \big[{-}\Et^T\Dat\Et\Fx(2:N_t+1, 1:N_x+1) - \Et^T\Dat\Ft + \Fx(2:N_t+1, 1:N_x+1)\B_x + \Et^T\A_4\big]\Ex^\dagger,
\end{align*}
we arrive again at
\begin{equation*}
\A\Uin + \Uin\B = \C.
\end{equation*}
In this case, once that $\Uin$ has been obtained, $\U$ follows from \eqref{e:Ugeneralboundary}.
	
\begin{remark} When only Dirichlet boundary conditions are considered, we have $d_a = d_b = 0$ in \eqref{e:mixuatubt}, and, hence, the boundary conditions in \eqref{e:mixuatubtmatrixnew} become $c_a\uu_{N_x+1} = \uu_a$ and $c_b\uu_1 = \uu_b$. Therefore, \eqref{e:Ex} and \eqref{e:Fx} get a much simpler form:
\begin{equation}
\label{e:ExFx}
\Ex =
\left(
\begin{array}{c|ccc|c}
	0 & 1 & & 0 & 0
	\\
	\vdots & & \ddots & & \vdots
	\\
	0 & 0 & & 1 & 0
\end{array}
\right), \qquad \Fx =
\left(
\begin{array}{c|ccc|c}
	\dfrac{u_b(t_0)}{c_b} & 0 & \ldots & 0 & \dfrac{u_a(t_0)}{c_a}
	\\
	\vdots & \vdots & \ddots & \vdots & \vdots
	\\
	\dfrac{u_b(t_{N_t})}{c_b} & 0 & \ldots & 0 & \dfrac{u_a(t_{N_t})}{c_a}
\end{array}
\right).
\end{equation}
Moreover, $\Ex\Ex^*$ and $(\Ex\Ex^*)^{-1}$ are the identity matrix of order $N_x - 1$, so $\Ex^\dagger = \Ex^*$ in that case.
\end{remark}

\subsection{Example with $x\in\mathbb R$}

\label{s:experiments}

In order to test our method for the case when $x$ is define on the whole real line, we have considered the following problem:
\begin{equation}
	\label{e:EDP1}
	\left\{
	\begin{aligned}
		& \cptt u(t,x) = u_{xx}(t,x) + 2xu_x(t, x) + 2u(t,x)
		\cr
		& \qquad\qquad\qquad - \frac{2^\alpha(\Gamma(1 - \alpha) - \Gamma(1 - \alpha, 2t))e^{2t - x^2}}{\Gamma(1 - \alpha)}, & & (t,x)\in[0, 1.2]\times\mathbb R,
		\cr
		& u(x, 0) = e^{-x^2},
	\end{aligned}
	\right.
\end{equation}
whose solution is $u(t,x) = e^{2t - x^2}$. In Listing~\ref{code:EDP1}, we offer its Matlab implementation. In order to generate the matrix $\Dat$, we invoke the function \verb|GenerateDa| defined in Listing~\ref{code:Da}; and we use the \verb|lyap| command to solve the resulting Sylvester equation. In order to measure the accuracy of the results, we define the error as $\|\operatorname{vec}(\U_{num}) - \operatorname{vec}(\U)\||_\infty$, where $\U_{num}$ denotes the numerical approximation obtained, $\U$ is the exact solution, and the $\operatorname{vec}$ operator piles up all the elements of a matrix in one single column.

\begin{lstlisting}[label=code:EDP1, language=Matlab, basicstyle=\footnotesize, caption = {Matlab program \texttt{EDP1.m}, that computes the numerical solution of \eqref{e:EDP1}}]
clear
tic
a=0.17; % $\alpha$
Nx=16; % $N_x$
Nt=2700; % $N_t$
tf=1.2; % $t_f$
t=tf*(0:Nt)'/Nt; % $t$
Dta=GenerateDa(Nt,tf,a); % $\mathbf D_t^\alpha$
b=1.4; % Scale factor $b$
[x,DD]=herdif(Nx,2,b);
x=x.'; % $x$
Dx=DD(:,:,1).'; % $\mathbf D_x$
Dx2=DD(:,:,2).'; % $\mathbf D_x^2$
[T,X]=ndgrid(t,x); % Two-dimensional grid
u=@(t,x)exp(2*t-x.^2); % $u(t, x)$
u0=@(x)exp(-x.^2); % $u_0(x)$
Et=[zeros(1,Nt);eye(Nt)]; % $\mathbf E_t$
Ft=[u0(x);zeros(Nt,Nx)]; % $\mathbf F_t$
A2=diag(2*x); % $\mathbf A_2$
A3=2*eye(Nx); % $\mathbf A_3$
A4=(2^a/gamma(1-a))*(gamma(1-a)-igamma(1-a,2*t)).*exp(2*T-X.^2); % $\mathbf A_4$
Bx=Dx2+Dx*A2+A3; % $\mathbf B_x$
A=Et'*Dta*Et; % $\mathbf A$
B=-Bx; % $\mathbf B$
C=-Et'*Dta*Ft+Et'*A4; % $\mathbf C$
Uinner=lyap(A,B,-C); % $\mathbf U_{inner}$
Unum=Et*Uinner+Ft; % $\mathbf U_{num}$
max(max(abs(Unum-u(T,X)))) % Error
toc % Elapsed time
\end{lstlisting}

In this example, we have chosen $\alpha = 0.17$, $N_t = 2700$ and $N_x = 16$. The code takes $1.84$ seconds to run, which includes also the cost of executing the \verb|igamma| function, and the error is only $1.6502\times10^{-10}$. Note that a value as small as $N_x = 16$ is enough to get spectrally accurate approximations in space, whereas $N_t$ much be much larger, because the error due to $\Dat$ behaves as $\mathcal O(N_t^{3-\alpha})$, as has been seen in Section~\ref{s:convergence}.

\subsubsection{Example with $x\in[x_a, x_b]$}
	
We have also considered a problem with $x$ belonging to a closed interval, with non-constant coefficients multiplying $u_{xx}$, $u_x$ and $u$, and taking time-dependent Robin boundary conditions:
\begin{equation}
	\label{e:EDP2}
	\left\{
	\begin{aligned}
		& \cptt u(t,x) = \frac{2^\alpha}{2.25}(1+x^2)u_{xx}(t,x) + \frac{2^\alpha}{1.5}x^2u_x(t,x)
			\cr
		& \qquad\qquad\qquad  - 2^{\alpha + 1}x^2u(t,x) - \frac{2^\alpha(\Gamma(1 - \alpha) - \Gamma(1 - \alpha, 2t))e^{2t + 1.5x}}{\Gamma(1 - \alpha)}, & & (t,x)\in[0, 1.2]\times[-1.1, 1.3],
		\cr
		& u(-1.1,t) + 2u_x(-1.1,t) = 4e^{2t-1.65},
		\cr
		& 3u(1.3, t) + 4u_x(1.3, t) = 9e^{2t+1.95},
		\cr
		& u(x, 0) = e^{1.5x},
	\end{aligned}
	\right.
\end{equation}
whose solution is $u(t,x) = e^{2t + 1.5x}$. In Listing~\ref{code:EDP2}, we offer its Matlab implementation. We invoke again the function \verb|GenerateDa| defined in Listing~\ref{code:Da}, and use the \verb|lyap| command to solve the Sylvester equation. Note that \verb|/Ex| in the definition of \verb|B| and \verb|C| is equivalent in Matlab to \verb|*pinv(Ex)|, where \verb|pinv(Ex)| is the Moore-Penrose inverse of \verb|Ex|, but \verb|/Ex| is more efficient.

\begin{lstlisting}[label=code:EDP2, language=Matlab, basicstyle=\footnotesize, caption = {Matlab program \texttt{EDP2.m}, that computes the numerical solution of \eqref{e:EDP2}}]
clear
tic
a=0.17; % $\alpha$
Nt=2700; % $N_t$
Nx=15; % $N_x$
xa=-1.1; % $x_a$
xb=1.3; % $x_b$
tf=1.2; % $t_f$
t=tf*(0:Nt)'/Nt; % $t$
Dta=GenerateDa(Nt,tf,a); % $\mathbf D_t^\alpha$
[Dx,x]=cheb(Nx);
x=xa+(xb-xa)*(x+1)'/2; % $x$
Dx=(2/(xb-xa))*Dx'; % $\mathbf D_x$
[T,X]=ndgrid(t,x); % Two-dimensional grid
u=@(t,x)exp(2*t+1.5*x); % $u(t, x)$
u0=@(x)exp(1.5*x); % $u_0(x)$
ca=1;da=2; % $c_a$ and $d_a$
ua=@(t)4*exp(2*t-1.65); % $u_a(t)$
cb=3;db=4; % $c_b$ and $d_b$
ub=@(t)9*exp(2*t+1.95); % $u_b(t)$
Et=[zeros(1,Nt);eye(Nt)]; % $\mathbf E_t$
Ft=[u0(x);zeros(Nt,Nx+1)]; % $\mathbf F_t$
c11=da*Dx(1,Nx+1);c12=ca+da*Dx(Nx+1,Nx+1); % $c_{11}$ and $c_{12}$
c21=cb+db*Dx(1,1);c22=db*Dx(Nx+1,1); % $c_{21}$ and $c_{22}$
Ex=[(-da*c22*Dx(2:Nx,Nx+1)+db*c12*Dx(2:Nx,1))/(c11*c22-c12*c21),eye(Nx-1),...
    (da*c21*Dx(2:Nx,Nx+1)-db*c11*Dx(2:Nx,1))/(c11*c22-c12*c21)]; % $\mathbf E_x$
Fx=[(c22*ua(t)-c12*ub(t))./(c11*c22-c12*c21),zeros(Nt+1,Nx-1),...
    (-c21*ua(t)+c11*ub(t))./(c11*c22-c12*c21)]; % $\mathbf F_x$
A1=diag((2^a/2.25)*(1+x.^2)); % $\mathbf A_1$
A2=diag((2^a/1.5)*x.^2); % $\mathbf A_2$
A3=diag(-2^(a+1)*x.^2); % $\mathbf A_3$
A4=-(2^a*igamma(1-a,2*t)/gamma(1-a)).*exp(2*T+1.5*X); % $\mathbf A_4$
Bx=Dx^2*A1+Dx*A2+A3; % $\mathbf B_x$
A=Et'*Dta*Et; % $\mathbf A$
B=-Ex*Bx/Ex; % $\mathbf B$
C=(-Et'*Dta*Et*Fx(2:end,:)-Et'*Dta*Ft+Fx(2:end,:)*Bx+Et'*A4)/Ex; % $\mathbf C$
Uinner=lyap(A,B,-C); % $\mathbf U_{inner}$
Unum=Et*Uinner*Ex+Et*Fx(2:end,:)+Ft; % $\mathbf U_{num}$
max(max(abs(Unum-u(T,X)))) % Error
toc % Elapsed time
\end{lstlisting}
	
In this example, we have taken $\alpha = 0.17$, $N_t = 2700$ and $N_x = 15$ (note again that a small value of $N_x$ is enough to get spectrally accurate approximations in space). The code takes $2.34$ seconds to run, and the error is only $1.8371\times10^{-10}$, which includes also the cost of executing the \verb|igamma| function.

\subsubsection{An example from the literature}

Finally, we have applied our method to the following equation taken from \cite[Example 3]{highorder3}:
\begin{equation}
	\label{e:EDP3}
	\left\{
	\begin{aligned}
		& \cptt u(t,x) = u_{xx}(t, x) - u_x(t, x) + \frac{720e^xt^{6 - \alpha}}{\Gamma(7 - \alpha)}, & & (t,x)\in[0, 1]^2,
		\cr
		& u(0, t) = t^6,
		\cr
		& u(1, t) = et^6,
		\cr
		& u(x, 0) = 0,
	\end{aligned}
	\right.
\end{equation}
whose solution is $u(t, x) = e^xt^6$. This example employs Dirichlet boundary conditions, and the initial data is the zero constant, so there are some simplifications in the code, which is offered in Listing~\ref{code:EDP3}. More precisely, $\Ex$ and $\Fx$ take the form in \eqref{e:ExFx}, with $c_a = c_b = 1$, $\Ex^\dagger = \Ex^T$, and $\Ft$ is an all-zero matrix, so it can be ignored.

\begin{lstlisting}[label=code:EDP3, language=Matlab, basicstyle=\footnotesize, caption = {Matlab program \texttt{EDP3.m}, that computes the numerical solution of \eqref{e:EDP3}}]
clear
tic
a=0.338; % $\alpha$
Nt=3500; % $N_t$
Nx=10; % $N_x$
tf=1; % $t_f$
t=tf*(0:Nt)'/Nt; % $t$
Dta=GenerateDa(Nt,tf,a); % $\mathbf D_t^\alpha$
[Dx,x]=cheb(Nx);
x=(x+1)'/2; % $x$
Dx=2*Dx'; % $\mathbf D_x$
[T,X]=ndgrid(t,x); % Two-dimensional grid
u=@(t,x)exp(x).*t.^6; % $u(t, x)$
ua=@(t)t.^6; % $u_a(t)$
ub=@(t)exp(1)*t.^6; % $u_b(t)$
Et=[zeros(1,Nt);eye(Nt)]; % $\mathbf E_t$
Ex=[zeros(Nx-1,1),eye(Nx-1),zeros(Nx-1,1)]; % $\mathbf E_x$
Fx=[ub(t),zeros(Nt+1,Nx-1),ua(t)]; % $\mathbf F_x$
A4=(720/gamma(7-a))*exp(X).*T.^(6-a); % $\mathbf A_4$
Bx=Dx^2-Dx; % $\mathbf B_x$
A=Et'*Dta*Et; % $\mathbf A$
B=-Ex*Bx*Ex.'; % $\mathbf B$
C=(-Et'*Dta*Et*Fx(2:end,:)+Fx(2:end,:)*Bx+Et'*A4)*Ex.'; % $\mathbf C$
Uinner=lyap(A,B,-C); % $\mathbf U_{inner}$
Unum=Et*Uinner*Ex+Et*Fx(2:end,:); % $\mathbf U_{num}$
max(max(abs(Unum-u(T,X)))) % Error
toc % Elapsed time
\end{lstlisting}

We have taken $N_t = 3500$, $N_x = 10$ and, as in \cite[Example 3]{highorder3}, we have considered three values of $\alpha$: $\alpha = 0.1$, whose error is $2.8880\times10^{-11}$; $\alpha = 0.2$, whose error is $1.6116\times10^{-10}$; and $\alpha = 0.338$, whose error is $7.2384\times10^{-10}$. In all the three cases, the elapsed time is approximately $2.50$ seconds. In general, ours errors are clearly smaller than those of \cite[Example 3]{highorder3}, which is due to the use of the Chebyshev differentiation matrices, instead of finite differences, to approximate the partial derivatives in space.

\section*{Funding}

Francisco de la Hoz was partially supported by the research group grant IT1615-22 funded by the Basque Government, and by the project PID2021-126813NB-I00 funded by MICIU/AEI/10.13039/501100011033 and by ``ERDF A way of making Europe''. Peru Muniain was partially supported by the research group grant IT1461-22 funded by the Basque Government, and by the project PID2022-139458NB-I00 funded by MICIU.

\end{document}